# Some isomorphic classes for noncanonical hypercomplex number systems of dimension 2

**Yakiv O. Kalinovsky,** Dr.Sc., Senior Researcher,
Institute for Information Recording National Academy of Science of Ukraine, Kyiv, Shpaka str. 2, 03113, Ukraine, E-mail: kalinovsky@i.ua

**Dmitry V. Lande,** Dr.Sc., Head of Department,
Institute for Information Recording National Academy of Science of Ukraine, Kyiv, Shpaka str. 2, 03113, Ukraine, E-mail: dwlande@gmail.com

**Yuliya E. Boyarinova,** PhD, Associate Professor,
National Technical University of Ukraine "KPI", Kyiv, Peremogy av. 37, 03056, Ukraine,
E-mail: ub@ua.fm

**Iana V. Khitsko,** Junior Researcher,
National Technical University of Ukraine "KPI", Kyiv, Peremogy av. 37, 03056, Ukraine,
E-mail: yannuary@yandex.ua

**Abstract**
Building of some isomorphic classes for noncanonical hypercomplex number systems o dimension 2 is described. In general case, such systems with specific constraints to structural constants can be isomorphic to complex, dual or double number system. Isomorphic transition between noncanonical hypercomplex number systems of the general form and diagonal form is built.

## Introduction

Structural data presenting forms are often used for calculations to increase the performance of different algorithms. One of these forms is hypercomplex data presenting form. More often canonical hypercomplex number systems are used, such as: complex, dual, double number systems, quasicomplex numbers, quaternions, octaves etc. However, for some tasks it is advisable to use non-canonical hypercomplex numbers, e.g. non-orthogonal coordinate systems for the analysis and synthesis of simple mechanisms **[1]**.

The most complete hypercomplex number system feature is a set of rules for basic elements multiplication. The system is called canonical if for each multiplication of every basic element pair the following condition is true:

$$E_i \cdot E_j = C_{ij} E_{ij}, \ C_{ij} \in \{-1, 0, 1\},$$

where $E_i, E_j, E_{ij}$ - are basic elements of hypercomplex number system. If at least one multiplication of basic element pair is a sum of two or more summands and/or has the coefficient out of $\{-1; 0; +1\}$, the HNS is called noncanonical **[1]**.

Multiplicity problems of noncanonical HNS are very important subjects for investigation, specifically enumeration problem **[2]** and building HNS isomorphism classes **[3]**.

## Isomorphism for noncanonical hypercomplex number system of dimenstion 2 in general

Let us consider the isomorphism of noncanonical HNS of 2$^{nd}$ dimenstion. **[4]**.

Further, hypercomplex number systems will be labeled as $\Gamma_i(\varepsilon, n)$, where $\varepsilon$ - system basic element, $n$ - system dimension.

In general, multiplication table of noncanonical commutative HNS of second dimension is as follows:

$$\Gamma_1(E,2) = \begin{array}{|c|c|} \hline a_{11}E_1 + b_{11}E_2 & a_{12}E_1 + b_{12}E_2 \\ \hline a_{12}E_1 + b_{12}E_2 & a_{22}E_1 + b_{22}E_2 \\ \hline \end{array}, \tag{1}$$

where $a_{11}, a_{12}, a_{22}, b_{11}, b_{12}, b_{22}$ - real numbers.

In order the multiplication table (1) to identify hypercomplex system, it must have the unit element $X = x_1 E_1 + x_2 E_2$ and hence the existence of non-trivial solution of hypercomplex equation:
$$M = X \cdot M, \tag{2}$$
where $M = m_1 E_1 + m_2 E_2 \in \Gamma_1$.

Then, taking into account (2.2) we are getting equations:

$$\begin{cases} (m_1 a_{11} + m_2 a_{12}) x_1 + (m_1 a_{12} + m_2 a_{22}) x_2 = m_1 \\ (m_1 b_{11} + m_2 b_{12}) x_1 + (m_1 b_{12} + m_2 b_{22}) x_2 = m_2 \end{cases}$$

Whence

$$x_1 = \frac{b_{12} m_1^2 + (b_{22} - a_{12}) m_1 m_2 + a_{22} m_2^2}{(a_{11} b_{12} - b_{11} a_{12}) m_1^2 + (a_{11} b_{22} - b_{11} a_{22}) m_1 m_2 + (a_{12} b_{22} - b_{12} a_{22}) m_2^2},$$

$$x_2 = \frac{-b_{11} m_1^2 + (a_{11} - a_{12}) m_1 m_2 + a_{12} m_2^2}{(a_{11} b_{12} - b_{11} a_{12}) m_1^2 + (a_{11} b_{22} - b_{11} a_{22}) m_1 m_2 + (a_{12} b_{22} - b_{12} a_{22}) m_2^2}.$$
(3)

It can be seen from (3), that the unit number depends on the number $M$ components, which should not be so. It is possible to get rid of such dependence if only quadratic forms coefficients in the nominators and denominators are proportional, i.e.:

$$\frac{b_{12}}{a_{11} b_{12} - b_{11} a_{12}} = \frac{b_{22} - a_{12}}{a_{11} b_{22} - b_{11} a_{22}} = \frac{a_{22}}{a_{12} b_{22} - b_{12} a_{22}},$$

$$\frac{-b_{11}}{a_{11} b_{12} - b_{12} a_{12}} = \frac{a_{11} - b_{12}}{a_{11} b_{22} - b_{11} a_{22}} = \frac{a_{12}}{a_{12} b_{22} - b_{12} a_{22}}.$$
(4)

Let us consider constants $a_{22}, b_{22}$ as free, since hypercomplex systems with such structural constants are already known **[4]**. Then we'll get the system of four equations with four unknown structural constants, which has such solutions:

1) $a_{11} \in R$, $b_{11} = 0$, $a_{12} = 0$, $b_{12} = a_{11}$, $a_{22}, b_{22} \in R$. (5)

This solution leads to a system, which is represented by the following multiplication table:

$$\Gamma_2(E,2) = \begin{array}{|c|c|} \hline a_{11}E_1 & a_{11}E_2 \\ \hline a_{11}E_2 & a_{22}E_1 + b_{22}E_2 \\ \hline \end{array}, \tag{6}$$

From unit (2) the following element of such system is got:
$$X = \frac{1}{a_{11}} E_1.$$

Using linear transformation of the $\Gamma_2(E,2)$ system basis, we can get its isomorphic system:

$$f_1 = \frac{1}{a_{11}} E_1; \quad f_2 = E_2.$$

$$\Gamma_3(f,2) = \begin{array}{|c|c|} \hline f_1 & f_2 \\ \hline f_2 & a_{11}a_{22} f_1 + b_{22} f_2 \\ \hline \end{array}, \qquad (7)$$

Then, depending on the value $a_{11}a_{22} + \frac{b_{22}^2}{4}$ sign, system (7) is isomorphic to the system of complex, double or dual numbers.

2) $b_{11} = 0$, $a_{12} = 0$, $a_{11} = b_{12} = b_{22}$, $a_{22} \in R \setminus \{0\}, b_{22} \in R$.

The following system is obtained

$$\Gamma_4(E,2) = \begin{array}{|c|c|} \hline b_{22}E_1 & b_{22}E_2 \\ \hline b_{22}E_2 & a_{22}E_1 + b_{22}E_2 \\ \hline \end{array}, \qquad (8)$$

with the unit element of:

$$X = \frac{1}{b_{22}} E_1.$$

Using linear transformations of the system (8) we can get its isomorphic system:

$$f_1 = \frac{1}{b_{22}} E_1; \quad f_2 = E_2.$$

$$\Gamma_5(f,2) = \begin{array}{|c|c|} \hline f_1 & f_2 \\ \hline f_2 & a_{22}b_{22} f_1 + b_{22} f_2 \\ \hline \end{array}, \qquad (9)$$

Depending on the given sign of the value

$$a_{22}b_{22} + \frac{b_{22}^2}{4},$$

system (9) is considered isomorphic to complex, double or dual number system.

### Isomorphism for noncanonical hypercomplex number system of dimension 2 and diagonal type

Let us consider some noncanonical hypercomplex number systems of the second dimension and diagonal type:

$$\Gamma_6(f,2) = \begin{array}{|c|c|} \hline \alpha_{11} f_1 + \beta_{11} f_2 & 0 \\ \hline 0 & \alpha_{22} f_1 + \beta_{22} f_2 \\ \hline \end{array}, \qquad (10)$$

Let the unit element for such system be $X = x_1 f_1 + x_2 f_2$, then it requires the existence of non-trivial real solution of equation (2), which turns into the linear equations system of:

$$\begin{cases} m_1 \alpha_{11} x_1 + m_2 \alpha_{22} x_2 = m_1 \\ m_1 \beta_{11} x_1 + m_2 \beta_{22} x_2 = m_2 \end{cases},$$

the solution of which is as follows:

$$x_1 = \frac{\beta_{22} m_1 m_2 + \alpha_{22} m_2^2}{(\alpha_{11} \beta_{22} - \beta_{11} \alpha_{22}) m_1 m_2},$$

$$x_2 = \frac{-\beta_{11} m_1^2 + a_{11} m_1 m_2}{(\alpha_{11} \beta_{22} - \beta_{11} \alpha_{22}) m_1 m_2}.$$

As in the previous case, we get rid of the $M$ number components dependence. It is required that the following condition is satisfied: $\beta_{11} = \alpha_{22} = 0$. In this case, system (10) takes the form of:

$$\Gamma_7(f,2) = \begin{array}{|c|c|} \hline \alpha_{11} f_1 & 0 \\ \hline 0 & \beta_{22} f_2 \\ \hline \end{array} \quad , \quad (11)$$

And the unit element takes the form of: $X = \dfrac{1}{\alpha_{11}} f_1 + \dfrac{1}{\beta_{22}} f_2$.

Let us define the conversion rules for the system $R \oplus R$ to system $\Gamma_7(f,2)$. Since both systems' multiplication tables are known, it is necessary to get the coefficients in isomorphism system:

$$R_1 * R_1 = R_1 = (y_{11} f_1 + y_{12} f_2)^2 = \alpha_{11} y_{11}^2 f_1 + \beta_{22} y_{12}^2 f_2 = y_{11} f_1 + y_{12} f_2$$

$$R_1 * R_2 = 0 = (y_{11} f_1 + y_{12} f_2) * (y_{21} f_1 + y_{22} f_2) = \alpha_{11} y_{11} y_{21} f_1 + \beta_{22} y_{12} y_{22} f_2 = 0,$$

$$R_2 * R_2 = R_2 = (y_{21} f_1 + y_{22} f_2)^2 = \alpha_{11} y_{21}^2 f_1 + \beta_{22} y_{22}^2 f_2 = y_{21} f_1 + y_{22} f_2$$

that gives us real equation system

$$\begin{cases} \alpha_{11} y_{11}^2 = y_{11} \\ \beta_{22} x_{12}^2 = y_{12} \\ \alpha_{11} y_{11} y_{21} = 0 \\ \beta_{22} y_{12} y_{22} = 0 \\ \alpha_{11} y_{21}^2 = y_{21} \\ \beta_{22} y_{22}^2 = y_{22} \end{cases},$$

which has two solutions:

1) $y_{11} = \dfrac{1}{\alpha_{11}}, y_{12} = 0, y_{21} = 0, y_{22} = \dfrac{1}{\beta_{22}};$

2) $y_{11} = 0, y_{12} = \dfrac{1}{\beta_{22}}, y_{21} = \dfrac{1}{\alpha_{11}}, y_{22} = 0.$

Both solutions give us a nonsingular linear transformation that converts system $R \oplus R$ to system $\Gamma_7(f,2)$. For the further analysis the second solution will be used.

$$R_1 = \dfrac{1}{\beta_{22}} f_2, R_2 = \dfrac{1}{\alpha_{11}} f_1. \qquad (12)$$

Les us build the isomorphic transition from system $\Gamma_1(E,2)$ to system $\Gamma_7(f,2)$ using intermediate transition to system $R \oplus R$. System $\Gamma_5(E,2)$ will be used as particular case of system $\Gamma_1(E,2)$. Then the transition operator of this kind of system to the system that is a direct sum of real numbers, will be as follows **[1]**:

$$\begin{cases} E_1 = R_2, \\ E_2 = -\dfrac{b_{22}}{2k} R_1 + \dfrac{1}{k} R_2, \end{cases} \qquad (13)$$

where $k^2 = a_{22}b_{22} + \dfrac{b_{22}^2}{4}$.

Taking into account (12)-(13) the transition operator from $\Gamma_5(E,2)$ to $\Gamma_7(f,2)$ is the following:

$$\begin{cases} E_1 = \dfrac{1}{\alpha_{22}} f_1, \\ E_2 = -\dfrac{b_{22}}{2k\beta_{22}} f_2 + \dfrac{1}{k\alpha_{22}} f_1. \end{cases}$$

## Conclusions

The article analysed some isomorphic classes for noncanonical hypercomplex number systems of dimension 2. It is shown that noncanonical hypercomplex number system of general form with the structure constants that correspond to the given constraints can be isomorphic to complex, double or dual numbers. The article also reveals the isomorphic transition from noncanonical diagonal system to a direct sum of real numbers, and as a result, the transition from noncanonical HNS of the general form to the noncanonical diagonal system.